\documentclass[11pt]{amsart}
\usepackage{amsmath,amsfonts}
\usepackage{amssymb}
\usepackage{amscd}
\usepackage{amsthm}
\usepackage{yhmath}
\usepackage{subfigure}
\usepackage[all]{xy}
\usepackage{color}

\numberwithin{equation}{section}

\newtheorem{theorem}[equation]{Theorem}

\newtheorem{proposition}[equation]{Proposition}

\newtheorem{lemma}[equation]{Lemma}

\newtheorem{corollary}[equation]{Corollary}

\theoremstyle{remark}
\newtheorem{remark}[equation]{Remark}

\theoremstyle{definition}
\newtheorem{definition}[equation]{Definition}

\newtheorem{question}[equation]{Question}

\def\XXint#1#2#3{{\setbox0=\hbox{$#1{#2#3}{\int}$}
	\vcenter{\hbox{$#2#3$}}\kern-.5\wd0}}

\newcommand{\Ad}{\operatorname{Ad}}

\newcommand{\To}{\longrightarrow}

\newcommand{\N}{\mathbb N}

\newcommand{\R}{\mathbb R}

\newcommand{\G}{{\mathbb G}}
\newcommand{\g}{\mathfrak{g}}

\newcommand{\norm}[1]{\left\Vert#1\right\Vert}

\begin{document}

\title[Isometries of Carnot groups and subFinsler 
homogeneous manifolds]{Isometries between open sets of Carnot groups and global isometries of subFinsler 
homogeneous manifolds}

\author{Enrico Le Donne}
\address{
Department of Mathematics and Statistics, University of Jyv\"askyl\"a, 40014 Jyv\"askyl\"a, Finland}
\email{enrico.ledonne@math.ethz.ch}

\author{Alessandro Ottazzi}
\address{
CIRM Fondazione Bruno Kessler, Via Sommarive 14, 38123 Trento, Italy}
\email{ottazzi@fbk.eu}

\renewcommand{\subjclassname}{%
 \textup{2010} Mathematics Subject Classification}

\date{\today}

\begin{abstract}
We show that isometries between open sets of Carnot groups are affine.
This result generalizes a result of Hamenst\"adt. Our proof does not rely on her proof. 
In addition, we  study  global isometries of general  homogeneous manifolds equipped with left-invariant subFinsler distances. We show that each isometry is determined by the  blow up at one point. 
For proving the results, we consider the action of  isometries  on the space of Killing vector fields. 
We make use of results by Capogna-Cowling and by Gleason-Montgomery-Zippin for obtaining smoothness of the isometric action.
\end{abstract}

\maketitle
\tableofcontents

%%%%%%%%%%%%%%%%%
%%%%%%%%%%%%%%%%%
\section{Introduction}

A fundamental problem in geometric analysis is the study of spaces 
that are isometrically homogeneous,
i.e., 
metric spaces on which the group of isometries acts transitively.
Such spaces have particular differentiable structures under the additional assumptions of being of finite dimension, locally compact, and the distance being intrinsic.
Indeed,
one can characterize such spaces as particular subFinsler manifolds, by using the theory of locally compact groups and methods from Lipschitz analysis on metric spaces, \cite{Gleason, Bochner_Montgomery46,
mz, b1, b2}.
Despite the fact that the group of global isometries of such manifolds is a Lie group and acts smoothly and by smooth maps,
the local isometries are still not completely understood.
Here, with the term `local isometry' we mean isometry between open subsets.

In this paper we give a complete description of the space of local isometries for those homogeneous spaces that also admit dilations.
These spaces, called Carnot groups, are particular nilpotent groups equipped with general left-invariant geodesic distances.

Our method of proof also shows that, as in Riemannian geometry, global isometries of homogeneous spaces are uniquely determined by their `first-order' expansion at a point.
Such a characterization for isometries was known already in some specific cases, e.g., for Riemannian manifolds, Banach spaces, and subRiemannian Carnot groups.

The study of isometries of distinguished Riemannian manifolds, such as homogeneous spaces, symmetric spaces, and Lie groups, has been a flourishing subject. 
 References for the regularity of isometries are the classical papers
\cite{Myers-Steenrod,Palais1}, see also \cite{Calabi-Hartman,Taylor}.
For the general theory of  transformations groups, we refer to  \cite{Palais, Kobayashi, 
Chen-Eberlein}.
Regarding the Finsler category, we mention the work  \cite{Deng-Hou}.
Banach spaces are classical and the space of isometries is well studied, see \cite{Fleming-Jamison, Benyamini-Lindenstrauss}.
There has been some effort in understanding isometries of subRiemannian manifolds, see \cite{Strichartz, Strichartz2, Hamenstadt, Kishimoto,Hladky}.
Regarding Carnot groups, in the deep paper \cite{Hamenstadt},
U.~Hamenst\"adt showed that  isometries are affine, in the case that the isometry is globally defined and that the distance is subRiemannian (and not just subFinsler).
We say that an isometry of a group (equipped with a left-invariant distance)
is
 affine  if it is the composition of a left translation with a group isomorphism.

We generalize 
Hamenst\"adt's result
to the setting 
of a subFinsler distance and isometries defined only on some open set.
%when 
%the distance is subFinsler and 
%the isometry is only defined
%on some open set.
 We need to point out that, first, to obtain such a local result, 
 one cannot use the same argument as in
 \cite{Hamenstadt} to deduce smoothness of the map.
 Actually, the issue of 
smoothness was a subtle point, which was clarified only later by I.~Kishimoto in \cite{Kishimoto}, for global isometries.
 Moreover, in  Hamenst\"adt's  strategy, 
 one needs to consider a blow down of the isometry, which requires the map to be   globally defined. 
 Hence, %for isometries between open sets, 
 we shall    provide a new method of proof.

%%%%%%%%%%%%%%%%%%%%%%%%%%%%%%%

\subsection{Statements of main results}
Let $G$ be a Lie group and $H$ be a closed subgroup.
Let $M=G/H$ be the homogeneous  manifold of left cosets. Hence, the group  $G$ acts transitively on $M$ on the left.
%Let $G$ be a Lie transformation group of an analytic manifold $M$.
Let $\Delta$ be a  $G$-invariant  subbundle of the tangent bundle $TM$.
We assume that $\Delta$ is bracket-generating and call it {\em horizontal bundle}.
%Consider a Lie group $G$  and a bracket-generating subspace $V$ of its Lie algebra $\g$.
Fix a norm on $\Delta_p$, for $p\in M$, and assume that it is  $G$-invariant. 
%Consider $V$ as a left-invariant subbundle of the tangent bundle.
Then the Carnot-Carath\'eodory distance between two points  of the manifold is 
the infimum of the lengths  of  curves
tangent to $\Delta$ and connecting the two points.
Since the length is measured using the norm, such a distance is also called {\em subFinsler}.
If $G=M$ and, setting $V_1:=\Delta_e$ and $V_{j+1}:=[V_1,V_j]$,
 we have the property  that  $${\rm Lie}(G)= V_{1}\oplus \cdots\oplus V_{s},$$
then the space is called {\em subFinsler Carnot group}. 
 %A subFinsler Carnot group is a   stratified Lie group equipped with such a Carnot-Carath\'eodory distance. 
% See Section \ref{prerequisites} for  more detailed definitions and properties of such spaces.
 %Finally, we fix the general notation and recall the prerequisites in Section \ref{prerequisites}.
 See Section \ref{prerequisites} for  more detailed definitions, notation, and properties of such spaces.

 % A subFinsler Carnot group is a nilpotent stratified Lie group equipped with a left-invariant subFinsler distance with the first stratum as horizontal distribution.

Our first theorem characterizes local isometries of a subFinsler Carnot group as affine maps.
\begin{theorem}\label{localisometries}
Let $(\G,d)$ be a subFinsler Carnot group.
Let $\Omega_1,\Omega_2\subset \G$ be two  open sets.
Let $f:\Omega_1 \to \Omega_2$ be an isometry.
Then  
there exists a left translation $\tau$ and a  group isomorphism $\phi $ of $\G$ 
such that $f$ is the restriction to $\Omega_1$ of $\tau\circ\phi$.
\end{theorem}
Note that in the statement above we require   the domain $\Omega_1$ to be open.
In Section \ref{afterthoughts}, we shall see  that such an assumption is necessary, unlike in the Euclidean case.
However, connectivity is not required.
%We don't require other topological assumptions such as connectivity.
%, without further assumptions. We shall see in Section \ref{afterthoughts} that if we do not suppose $\Omega_1$ to be open, then the conclusion of the theorem may fail to hold.

With a similar method as for the proof of the first theorem, we show that global isometries of a subFinsler  homogeneous space are characterized by their 
value at one point and the differential at this point.
In \cite{Strichartz, Strichartz2}, the same conclusions were obtained  for smooth isometries of particular subRiemannian manifolds.

%We denote by $e$ the identity of $G$. % and by $H_eG$ the horizontal space at the identity.
\begin{theorem}\label{globalisometries}
Let $M=(G/H,d)$ be a connected homogeneous
manifold equipped with a $G$-invariant subFinsler distance with horizontal subbundle $\Delta$. 
Let $f:M\to M$ be an isometry.
Then $f$ is an analytic map. % and it is completely determined by its blow up at the identity $e$. Namely, 
Moreover, if $h:M\to M$ is another isometry with the properties that
$f(p)=h(p)$ and $(df)_p|_{\Delta_p}=(dh)_p|_{\Delta_p}$, for some $p\in M$, then $f=h$.
\end{theorem}

Even under the additional assumption that the manifold in Theorem  \ref{globalisometries}  is actually a Lie group, one cannot have the same conclusions as in Theorem 
\ref{localisometries}.
%with the assumptions of Theorem  \ref{globalisometries}.
 Indeed, for general Lie groups, it is not true that isometries  are necessarily affine maps.
  See the discussion on the inversion map on $\mathbb S^3$ in Section \ref{afterthoughts}.
We should also point out that,  in general, local isometries of homogeneous spaces do not extend to global isometries, e.g., in the case of the flat cylinder (see Section \ref{afterthoughts}).  
  
By \cite{Margulis-Mostow}, every quasiconformal mapping 
between two Carnot-Carath\'eodory spaces
admits a blow-up map at almost every point,
that is an isomorphism between two Carnot groups. 
% which is and isomorphism between the blow-up spaces. 
%The latter turn out to be Carnot groups \cite{Mitchell}. 
%If the manifold is a Lie group and the distance is left-invariant, then the blow-up space is the same at every point. 
%In particular, if we start with a quasiconformal transformation  between domains in a subFinsler Carnot group, then the blow-up map is also called the Pansu differential \cite{Pansu}. 
If the two Carnot-Carath\'eodory spaces are Carnot groups, such a result was proven in \cite{Pansu}. 
For this reason,
the
blow-up map is also called   Pansu differential.
Whenever $f$ is an isometry of $M$, as in Theorem \ref{globalisometries}, then $f$ is smooth and $df(p)|_{\Delta_p}$ coincides with the differential of the blow up at $p$, restricted to $\Delta_p$. 
In other words, Theorem  \ref{globalisometries} states that every isometry $f$ of $M$ is determined by $f(p)$ and by its blow-up at $p$.

%Pansu differenziabile? Blow up?

We complete the introduction with a problem.
\begin{question}   
Given   a connected and simply connected Lie group $(G,d)$ equipped with a left-invariant subFinsler distance, let $\G$ be the tangent cone of $G$ at $e$. For which group automorphism $\phi$ of $\G$  there exists an isometry of $G$ that has $\phi$ as  blow up at the identity?
\end{question}   
There are obvious necessary conditions. Namely, the differential of $\phi$ should preserves the strata and be an isometry when restricted to the first stratum. However, these conditions are not sufficient, see the discussion on the Riemannian Heisenberg group in Section \ref{afterthoughts}.

\subsection{Structure of the paper}

The proof of our results is divided into several steps. 
The overall strategy of the two theorems is similar. However, in few instances we need different approaches.
Regarding  global isometries of subFinsler homogeneous manifolds, we show their smoothness in Section \ref{globalsmoothness}, relying on the classical solutions of the Hilbert $5^{th}$   problem. This method works only for mappings that are globally defined. For the local isometries of Theorem \ref{localisometries} we shall rely on a regularity result of $1$-quasiconformal mappings in \cite{Capogna-Cowling}, that we reformulate in Theorem \ref{capcow}. 

Once we know that isometries are smooth, in Section \ref{killing} we consider the action of their differentials on vector fields that generate flows of isometries (Killing vector fields).  
For local isometries, we need an extension result for Killing vector fields
 that is proven in Section \ref{tanaka} and that relies on a method developed in \cite{Tanaka}. 
 In Section \ref{globalthm}  we prove Theorem \ref{globalisometries}.
 In Section  \ref{localthm}
  we provide  a weak version of Theorem 
\ref{localisometries}.  Namely, 
in Theorem \ref{local-connected},
we assume that $\Omega_1$ is connected and that the distance is subRiemannian. 
In Section \ref{endproof}
 we complete the proof of Theorem 
\ref{localisometries}.    We then devote Section \ref{afterthoughts} to a number of final remarks.

%Finally, we fix the general notation and recall the prerequisites in Section \ref{prerequisites}.

\subsubsection{Acknowledgements}
Both authors would like to thank   Universit\'{e} Paris Sud,  Orsay, % for the excellent working environment, %for perfect working conditions, 
where part of this research was conducted.
This paper has benefited  from discussions with  E. Breuillard and P. Pansu.
Special thanks go to them. 

%%%%%%%%%%%%%%%%%%%%%%%%%%%%%%%%%%%%

%%%%%%%%%%%%%%%%%%%%%%%%%%%%%%%%%%%%%%%%%%%%%%%%
%%%%%%%%%%%%%%%%%%%%%%%%%%%%%%%%%%%%%%%%%%%%%%%%%

\section{Preliminaries}\label{prerequisites}
%%%%%%%%%%%%%%%%%%%%%%%%%%%%

\subsection{General notation}\label{General_notation}
Let $G$ be a Lie group. Denote by $\g$ or by Lie$(G)$ the Lie algebra of $G$ whose elements are tangent vectors at the identity $e$ of $G$.
For $Y\in \g$, we denote by $\tilde Y$   the left-invariant vector field that coincides with $Y$ in $e$. So $[X,Y]:=[\tilde X,\tilde Y]_e$.

Let $H$ be a closed subgroup of $G$. Hence, the space $G/H$ of left cosets $gH$, with $g\in G$, has a natural structure of analytic manifold, see  \cite[page 123]{Helgason}.  
The group $G$ is a Lie transformation group of   $M =G/H$. Namely, every element $g\in G$ acts by left translations on $M$, i.e., induces the diffeomorphism
\begin{equation}\label{azione}
g'H\mapsto g\cdot(g'H):=gg'H.
\end{equation}

For $Y\in \g$, we denote by  $Y^\dag$  the   vector field of $M$ whose flow $\Phi_{Y^\dag}^t$ at time $t$ is
\begin{equation}\label{left-group}
\Phi_{Y^\dag}^t(p)=\exp(tY)\cdot p, \qquad \forall p\in M.
\end{equation}
%For $Y\in \g$, we denote by $\tilde Y$ and $Y^\dag$  the left- and the right-invariant vector field that coincides%with $Y$ at the identity $e$, respectively.
It is well known (see \cite[Theorem 3.4]{Helgason})
 that, for $X,Y\in\g$, we have 
 \begin{equation}\label{bracket_dag}
[X,Y]^\dag=-[X^\dag,Y^\dag].
\end{equation}

We shall fix a $G$-invariant subbundle $\Delta$ of the tangent bundle $TM$ of $M$. The choice of such a subbundle can be seen in the following way.
In the homogeneous manifold $M$ we denote by $o$ the coset $H$ and call it the {\em origin} of $M$. Notice that the action of $H$ on $M$ fixes the origin.
There is a one-to-one correspondence between $H$-invariant subspaces $\Delta_o$ in  $T_o(M)$ and $\Ad_H$-invariant subspaces $V$ in Lie$(G)$ that contains Lie$(H)$.
We choose such a  subspace $\Delta_o$ in  $T_o(G/H)$, and therefore, such a $V\subseteq \g$.
Then, for all $gH\in G/H$, the subbundle $\Delta$ is defined as
$$\Delta_{gH} :=   g_*  \Delta_o  
, $$
where $g_*$ is the differential of the map in \eqref{azione}. The subbundle is well defined, i.e., the definition does not depend on the representative in $gH$,  exactly because $\Delta_o$ is $H$-invariant. 

%{\color{red} 
%Let $\tilde \Delta_g:=(d\pi_g)^{-1} \Delta_o$.
%Hence $V=\tilde \Delta_e$.
%Claim $\tilde \Delta$ is $\Ad_H$-invariant.
%
%D'altra parte, se $\tilde \Delta$ e' $G$-invariant e definitamo $d\pi_e \tilde \Delta_e=: \Delta_o$, allora $\Delta_o$ si estende a un $G$-invariant subbundle $\Delta$ se (e solo se) $\tilde \Delta$ is $\Ad_H$-invariant.
%
%Quindi: mi sa che $V$ deve essere scelto anche $\Ad_H$-invariant.}

If the  subspace $V\subset \g$ associated to $\Delta_0$   has the property that $\g$ is the smallest Lie subalgebra of $\g$ containing $V$, then $V$ (or, equivalently, $\Delta$) is said to be {\em bracket-generating}.

%Here and hereafter,  
%we use both notations
% $df$ or  $f_*$
% to denote 
%  the   differential of a differentiable map $f$, viewed as 
% push-forward operator on vectors or on  vector fields respectively.
%Denote by $L_p: G\to G$ the left translation by  $p\in G$. Define
%$$\tilde V _p :=   (dL_p)_e  V  
%.  $$
%Hence 
%$\tilde V$ is  the left-invariant  subbundle of the tangent bundle $TG$ such that $\tilde V_e =V$.

We shall fix a $G$-invariant norm on $\Delta$.  The choice of such a norm can be seen in the following way.
Fix a seminorm on $V$ that is $\Ad_H$-invariant and for which the kernel is Lie$(H)$.
%Fix a norm    on $V$ that is $\Ad_H$-invariant.
The projection from $G$ to $M$ gives an $H$-invariant norm  $\norm{\cdot}$ on $\Delta_o$.
Hence, we have an induced 
$G$-invariant norm on $\Delta$  by 
$$\norm{v}=\norm{({g}^{-1})_*   v} , \qquad \forall v\in \Delta_{gH}.$$
Since the initial norm is $\Ad_H$-invariant, it follows that the above equation is independent on the choice of representative in $ gH$.

An absolutely continuous
curve $\gamma:[0,1]\to M$   is said to be   {\em horizontal} (with respect to  $\Delta$) if the derivative
 $\dot\gamma(t)$ belongs to $\Delta$, for almost every $t\in [0,1]$.\
 Each horizontal curve $\gamma$
 has an associated  length
 defined as  %follows. 
 %Let $\sigma(t)\in G$ such that $\gamma(t)=\sigma(t)H$. The value $\sigma(t)$ might not be continuous in $t$. One define the length of $\gamma$ as
$$
    L(\gamma) := 
    \int_0^1 \norm{\dot\gamma(t)}\,\mathrm d t
    .
$$
\begin{definition}[SubFinsler   homogeneous manifolds]\label{SF_manifold}
Let $M=G/H$ be a homogeneous space formed by  a Lie group $G$ modulo a closed subgroup $H$.
We are given a bracket-generating
$G$-invariant subbundle  $\Delta \subseteq TM$ and a $G$-invariant norm $\norm{\cdot}$ on $\Delta$.
Equivalently, we are given 
an $\Ad_H$-invariant and   bracket-generating subspace $V\subseteq {\rm Lie}(G)$, with $V\supseteq {\rm Lie}(H)$, and an $\Ad_H$-invariant seminorm  $\norm{\cdot}$ on $V$ whose kernel is Lie$(H)$.
 The  {\em subFinsler  distance} (also known as Finsler Carnot-Carath\'eodory distance) between two points $p,q\in M$ is defined as
 \begin{equation}\label{dist_sF}
d(p,q):=\inf\{L(\gamma)\;|\;\gamma\text{ horizontal and } \gamma(0)= p, \gamma(1)= q\}.
\end{equation}
We call the pair $(M,d)$ a 
subFinsler   homogeneous manifold.
\end{definition}
By Chow's Theorem \cite[Chapter 2]{Montgomery}, the topology of  $(M,d)$ is the  topology of $M$ as manifold.
Notice that, by construction,  the above subFinsler distance
is left-invariant, i.e., every left translation \eqref{azione} is an isometry of $(M,d)$.

By the work of V.~N.~Berestovski{\u\i}, %and the earlier work of Montgomery and Zippin, 
we know that the above-defined subFinsler   homogeneous manifolds are the only geodesic spaces that are isometrically homogeneous, are locally compact, have finite topological dimension, and whose distance is a geodesic distance. 
Such a  result is based on Montgomery-Zippin's characterization of Lie groups, see Theorem \ref{Mont-Zip2}.

\begin{theorem}[Consequence of \cite{mz}, \cite{b2}, and \cite{Mitchell}] \label{Bere} 
Let $X$ be a locally compact and finite-dimensional topological space.
Assume that $X$ is equipped with an intrinsic distance $d$ such that
its isometry group ${\rm Iso}(X,d)$ acts transitively 
on $X$.
Then $(X,d)$ is isometric to a subFinsler   homogeneous manifold.

If, moreover, the space $(X,d)$ admits a non-trivial dilation, i.e., there exists $\lambda>1$ such that
$(X,\lambda d)$ is isometric to $(X,d)$, then $(X,d)$ is a subFinsler Carnot groups (see definition below).
\end{theorem}

 By the above result, subFinsler Carnot groups are special cases of subFinsler     homogeneous manifold.  We refer to \cite[page 38]{Montgomery} for the easy proof that such spaces admit dilations, for all $\lambda >1$.
 
 %Given a subspace $V$ of the Lie algebra of a group, we denote by $\tilde V$    the left-invariant  subbundle of the tangent bundle   such that $\tilde V_e =V$.

\begin{definition}[SubFinsler Carnot groups]
Given a subspace $V_1$ of the Lie algebra of a Lie group $G$, define by recurrence the subspaces $V_j$ as
$$V_j:=[V_1,V_{j-1}], \qquad \forall j>1.$$
If one has that 
$$\g= V_{1}\oplus \cdots\oplus V_{k}\oplus \cdots, $$
then $G$ is said to be a (nilpotent)  {\em stratified group} and
 $V_1$ is called the {\em first stratum} (of the stratification $\{V_j\}$).
 If $d$ is the 
Carnot-Carath\'eodory  distance
associated to $G$, $V_1$, and some norm $\norm{\cdot}$ on $V_1$, then the pair 
$(G,d)$ is called  
{\em subFinsler Carnot group}, or simply Carnot group.
\end{definition}

If the norm in Definition \ref{SF_manifold} comes from a scalar product, then the associated distance is called {\it Carnot-Carath\'eodory} or {\it subRiemannian}. 
If this is the case for a subFinsler Carnot group, then we call it {\em subRiemannian Carnot group}. 
We shall use the notation $\G$, rather than $G$,  to emphasize that we are dealing with a  Carnot group, rather than a general Lie group.

One can show that a curve in a subFinsler     manifold has finite length  if and only if it is a horizontal curve, up to reparametrization.
Consequently in the case that an isometry $f$ of a subFinsler   homogeneous manifold   %$(G,d)$ 
is $C^1$, then it is a {\it contact map}, i.e., its differential preserves the subbundle. Namely,
$$df_p (\Delta_p) \subseteq \Delta_{f(p)}, \qquad \forall p\in G.$$
Here and hereafter,  
we use both notations
 $df$ or  $f_*$
 to denote 
  the   differential of a differentiable map $f$, viewed as 
 push-forward operator on vectors or on  vector fields.

We shall need to show that isometries between open sets of a subRiemannian Carnot group are analytic maps. Such a regularity result 
follows from the work of Capogna and Cowling on $1$-quasiconformal mappings.
We state here a weaker form of their result (\cite[Theorem 1.1]{Capogna-Cowling}).

\begin{theorem}[Consequence of \cite{Capogna-Cowling}]\label{capcow}
Let $\Omega\subseteq \G$ be an open set of a subRiemannian Carnot group $\G$. Let $f:\Omega\to \G$ be a biLipschitz embedding.
\begin{itemize}
\item[(i)] If $f$ is an isometry, then $f$ is analytic.
\item[(ii)] If for a.e. $p\in \Omega$, the blow-up of $f$ at $p$ is an isometry of $\G$, then $f$ is analytic.

\end{itemize}
\end{theorem}

%%%%%%%%%%%%%%%%%%%%%%%%%%%%
\subsection{Smoothness of the isometric action}\label{globalsmoothness}

By definition, an isometry is a map that preserves the distance. Hence, there is no a priori  smoothness assumption.
In this section, we shall explain why, in the case of subFinsler homogeneous manifolds, one has in fact that (global) isometries are smooth maps forming a Lie group and the action is smooth. Smoothness of the local action of
 isometries between open sets of Carnot groups  will follow from a different reasoning, see Corollary \ref{Hilbert}.

The smoothness of global maps is a consequence of the work of Gleason  \cite{Gleason}  and Montgomery and Zippin \cite{Montgomery_Zippin52, mz}.
In particular, one has the following general result.
\begin{theorem}[Gleason-Montgomery-Zippin] \label{Mont-Zip2} 
If a second countable and  locally compact group $H$ acts by isometries, continuously, effectively, and transitively on a locally compact, locally connected, and finite-dimensional metric space $X$ then $H$ is a Lie group and $X$ is a differentiable manifold.
\end{theorem}

To obtain the regularity of the action in the group parameters we use  the following theorem, which is a generalization of Bochner-Montgomery's result \cite{Bochner_Montgomery45}.
\begin{theorem}[{\cite[page 213]{mz}}]\label{teoA}
Let $H\times M\to M$, $(h,x)\mapsto h(x)$, be a (continuous) action of a Lie group $H$ on an analytic manifold $M$.
Assume that, for all $h\in H$,  the map $x\mapsto h(x)$ is analytic.
Then $h(x)$ is analytic in $h$ and $x$ simultaneously. 
\end{theorem}

For studying local isometries, we will make use of another result of Montgomery, see \cite{Montgomery_45_1,Montgomery_45_2}. We obtain a Lie group structure for {\em compact} groups acting by analytic maps. The result holds more generally, see \cite[page 208, Theorem 2]{mz}, but we only need the following weaker result.
\begin{theorem}[{\cite[Theorem 3]{Montgomery_45_1}}]\label{teoB}
If $H$ is a compact effective group acting on a connected analytic manifold $M$ and if each transformation of $H$ is analytic then $H$ does not contain arbitrarily small subgroups other than the identity; or, in other words, $H$ considered in itself is a Lie group.
\end{theorem}

From Theorem \ref{Mont-Zip2},  we deduce the following consequence, which was partially observed in \cite{Kishimoto} as well.
\begin{corollary}[Consequence of Hilbert 5$^{th}$ theory]\label{compactness}
Let $G$ be a  Lie group acting  on an analytic manifold $M$.
Assume that the action is transitive and analytic.
Let $d$ be a  $G$-invariant distance on $M$, inducing the manifold topology.
Then the isometry group ${\rm Iso} (M)$ is a Lie group, 
the action
 \begin{eqnarray}\label{Iso_azione}
 {\rm Iso} (M)\times M&\to& M\\
 (f,p)&\mapsto& f(p) \nonumber
 \end{eqnarray}
 is analytic, and the space 
 $$
{\rm Iso}_o (M)=\{f \in {\rm Iso}(G) \mid f(o)=o\}
$$
is a compact Lie group.
\end{corollary}

\proof
By Ascoli-Arzel\`a Theorem  we have that ${\rm Iso} (M)$ is  locally compact and $
{\rm Iso}_o (M)$ is compact (both equipped with the compact open topology).
Obviously they both are  groups with the composition as multiplication.
Furthermore, 
since $G$ acts transitively on $M$, so does ${\rm Iso} (M)$. 
Therefore, by Theorem \ref{Mont-Zip2} it follows that ${\rm Iso} (M)$ is a Lie group.
Being a compact subgroup, $
{\rm Iso}_o (M)$ is a Lie group as well.

% =========================================
 
For the proof that the action of ${\rm Iso} (M)$ on $M$ is analytic, 
% we should show that $M$ 
 let us  explicit  the analytic structures considered.
 The group $G$ and the manifold $M$ are given with their analytic structures, which we denote
 $\omega_G$ and $\omega_M$, respectively.
 Hence, by assumption, the action
 \begin{equation}\label{prima_mappa}
 (G,\omega_G)\times (M,\omega_M) \To (M,\omega_M),
 \end{equation}
given by \eqref{azione}, is analytic.
The group ${\rm Iso} (M)$ 
has an analytic structure $ \omega_I$ of Lie group and, since it 
is acting transitively (and continuously) on $M$,
there exists an  analytic structure $\tilde\omega_M$ on $M$ such that the map
 \begin{equation}\label{seconda_mappa}
 (I,\omega_I)\times (M,\tilde\omega_M) \To (M,\tilde\omega_M),
 \end{equation}
given by \eqref{Iso_azione}, is analytic, see \cite[page 123]{Helgason}. Every element of $G$ induces an isometry. Hence, we have a map 
  \begin{equation}\label{quinta_mappa}
\iota: (G,\omega_G)\To (I,\omega_I),
 \end{equation}
 induced by \eqref{azione}.
The map $\iota$ is a continuous homomorphism. 
 By \cite[Theorem 2.6]{Helgason} we have that  $\iota$ is in fact analytic.
By composition of \eqref{seconda_mappa} and  \eqref{quinta_mappa}, we have that 
 \begin{equation}\label{quarta_mappa}
(G,\omega_G)\times (M,\tilde\omega_M) \To (M,\tilde\omega_M),
 \end{equation}
again given by \eqref{azione},  is analytic.
 By \cite[Theorem 4.2]{Helgason} 
there is a unique analytic structure on $M$ for which the action given by  \eqref{azione} is analytic.
Therefore, we conclude that $\omega_M=\tilde \omega_M$.
Hence, the map \eqref{Iso_azione} is analytic when $M$ is given the initial analytic structure $\omega_M$.
\qed

From  Theorem \ref{teoB} we deduce the following consequence, which will be important in the study of local isometries of Carnot groups.
\begin{corollary}[Consequence of Hilbert 5$^{th}$ theory]\label{Hilbert}
Let $B$ be a closed ball in a Carnot group centered at the identity $e$.
Let $\mathcal I :=\{f:B\to B  \,\,{\rm isometry} \mid f(e)=e \}.$
Then \begin{itemize}
\item[(i)]   $\mathcal I$ is compact;
\item[(ii)] $\mathcal I$ is a Lie group;
\item[(iii)] for all $x\in B$ the map
 \begin{eqnarray*}
 \mathcal I &\to& B\\
 f&\mapsto& f(x)
 \end{eqnarray*}
 is analytic.
\end{itemize}
\end{corollary}

\proof
%(i) Ascoli-Arzela
%(ii) by \cite{Capogna-Cowling} and thm 4.39
%(iii) thm 4.51
%
%OR use Neumann, Smith ?
Regarding part (i), since $\mathcal I$ is a family of equicontinuous and equibounded maps, by
 Ascoli-Arzel\`a Theorem, the family is precompact. Since $\mathcal I$ is obviously closed, then it is compact.
 
 Part (ii) will follow from Theorem \ref{teoB}. Indeed, the set $\mathcal I$ is closed under composition. Hence it is a compact group acting on $B$.
 By  Theorem \ref{capcow}, each element of $\mathcal I$ is an analytic map.
Then Theorem \ref{teoB} implies that 
$\mathcal I$  is a Lie group.

Part (iii) follows immediately from Theorem \ref{teoA}.  
\qed

\subsection{The action of an isometry on Killing vector fields}\label{killing}

%For $Y\in \g$, we denote by $\tilde Y$ and $Y^\dag$ respectively the left and the right-invariant vector field  that %coincides with $Y$ at the identity $e$.

%Moreover, we write $\tilde V_1$ to indicate the sub-bundle of the tangent bundle that is obtained by left translating %$V_1$.

 %For a map $f$ between subsets of $G$, we denote by $df$ its differential whenever this exists. We shall use the %notation $f_*$ to indicate the push forward of vector fields.
%\begin{definition}
%Given open and connected subsets $U,V\subseteq G$, a map $f: U\to V$ is an {\it isometry} if $d(f(p),f(q))=d(p,q)$, for %every $p,q\in U$.
%\end{definition}
%It is worth to notice that if an isometry $f$ is differentiable, then in particular $df (\tilde V_1) \subseteq \tilde V_1$, in %other words $f$ is {\it contact}.

Let now $M=(G/H,d)$ be a subFinsler homogeneous space with horizontal bundle $\Delta$.
In this section we define a filtration of the space of (global) vector fields that are infinitesimal generators of isometries of $M$. 
The properties of this space that are shown here are crucial to prove both  our results.

\begin{definition}[Killing vector fields $\mathcal K$]
A vector field $Z$ on $M$ is said to be a {\em Killing vector field} if there exists $t_0>0$ such that, for all $t\in [0,t_0]$, the flow $\Phi_Z^t$ at time $t$ is an isometry.
We denote by $\mathcal K$ the collection of all  Killing vector fields.
\end{definition}
One can show that $\mathcal K$ is closed under sum and Lie bracket. 
Notice that if $Y\in V$, then the vector field $Y^\dag$ is a Killing vector field. Indeed, the flows of these vector fields are one parameter groups of  the left action of $G$ (see \eqref{left-group}).
The space $\mathcal K$ is the Lie algebra of the group of (global) isometries of $M$. Such a group, by Corollary \ref{compactness}, is a Lie group. Hence, $\mathcal K$ is a finite dimensional Lie algebra.
We recall that isometries that are smooth are in particular contact maps. Therefore, given a Killing vector field $Z$ with flow $\Phi_Z^t$, the fact that $\Phi_Z^t$ is contact implies that
\begin{equation}\label{contactvectorfield}
[Z,\Gamma ( \Delta)] \subset \Gamma(\Delta),
\end{equation}
where $\Gamma(\Delta)$ denotes the space of smooth sections of the subbundle $\Delta$.
Indeed, if $W\in\Gamma(\Delta)$, we have
$$[W, Z]=  -\mathcal L_Z (W) =
\frac{d}{dt} (\Phi_Z^t)_*  (W)|_{t=0} \in \Gamma(\Delta) ,
$$
where $\mathcal L_Z (W)$ denotes the Lie derivative along $Z$ of $W$.

As in Section \ref{General_notation}, we denote by $o$the origin in $M$.
We write $\mathcal K_{-1}:=\{Z\in \mathcal K\mid Z_o\in \Delta_o\}.$ Moreover,
  we denote $\mathcal K_0:=\{Z\in \mathcal K \mid Z_o =0\}$ and inductively, for every $j\geq 1$, we define
$$
\mathcal K_j := \{Z\in \mathcal K_0 \mid [Z,Y^\dag]\in \mathcal K_{j-1} \,\forall Y\in V\}.
$$
Notice that $\mathcal K_{i+1} \subseteq \mathcal K_i$ for every $i\geq -1$. Moreover, each $\mathcal K_j$ is a vector space.

If $Z$ is a Killing vector field, we choose $Y\in \g$ such that $(Y^\dag)_o = Z_o$ and therefore we decompose $Z$ as 
\begin{equation}\label{decomposition}
Z= Y^\dag + Z^\prime, 
\end{equation}
with $Z^\prime\in \mathcal K_0$. We obtain the decomposition  $\mathcal K= \g^\dag + \mathcal K_0$, which is not necessarily a direct sum.
We show the following property.

\begin{lemma}\label{gradationofkilling}
$[\mathcal K_0, \mathcal K_j]\subseteq \mathcal K_{j}, \,\,\forall j \geq 0$.
\end{lemma}
\proof
We prove the lemma using induction. If $j=0$, the statement is true, because the bracket of two vector fields that vanish at $0$ vanishes at $0$ as well. 

Suppose now that the claim is verified for all indexes from $0$ to $j-1$ and pick $Z\in \mathcal K_0$ and $Z^\prime\in \mathcal K_j$.  Clearly $[Z,Z^\prime]_o=0$, i.e., $[Z,Z^\prime]\in \mathcal K_0$. It remains to prove that $[[Z,Z^\prime],Y^\dag]\in \mathcal K_{j-1}$, for every $Y\in V$.
First notice that the Jacobi identity yields
\begin{equation}\label{justonesum}
[[Z,Z^\prime],Y^\dag]= -[[Z,Y^\dag],Z^\prime] + [[Z^\prime,Y^\dag],Z].
\end{equation}
Recall the identification between  $V$ and $\Delta_o$ given in Section \ref{prerequisites} and let $Y^\prime\in\Gamma(\Delta)$ be such that $Y^\prime_o = Y$.
Then $[Z,Y^\dag]= [Z,Y^\prime]+ [Z,Y^\dag -  Y^\prime]$. Therefore, it follows that $[Z,Y^\dag]_o =[Z, Y^\prime]_o$. We recall that $Z$ satisfies \eqref{contactvectorfield}, so that in particular $[Z,Y^\dag]_o \in \Delta_o$. Thus \eqref{decomposition} gives $[Z,Y^\dag]=W^\dag + W^\prime$, for some $W\in V$ and $W^\prime\in \mathcal K_0$. Since $\mathcal K_j\subseteq \mathcal K_{j-1}$, by induction we conclude that
\begin{equation}\label{righthand}
[[Z,Y^\dag],Z^\prime] = [W^\dag,Z^\prime]   + [W^\prime,Z^\prime] \in \mathcal K_{j-1}.
\end{equation}
Again by the induction hypothesis we have that $[[Z^\prime,Y^\dag],Z] \in \mathcal K_{j-1}$, which, together with \eqref{righthand} and \eqref{justonesum}, finishes the proof.
\qed

\begin{lemma}\label{gradinggeneral}
Let $f$ be an isometry between open and connected subsets of  $M$. 
We assume that $f$ is smooth and  that $f_*Z$ uniquely extends to a Killing vector field for every $Z\in\mathcal K$, for which we shall abuse the notation $f_*Z$.
Moreover, we suppose that $df_o|_{\Delta_o}$ is the identity, and that $f(o)=o$. Then
\begin{itemize}
\item[(i)]   if $Z\in \mathcal K_{-1}$, then $f_*Z\in Z + \mathcal K_0$;\\
\item[(ii)]  if $Z\in \mathcal K_j$ with $j\geq -1$, then $f_* Z\in Z +\mathcal K_{j+1}$.

\end{itemize}
\end{lemma}
\proof
The idea of  the proof is the following. The first part is a consequence of the fact that $df_{o}$ is the identity on $\Delta_o$. The second part will follow by induction and Lemma \ref{gradationofkilling}.

Regarding the proof of (i), 
we first note that  by hypothesis $f_*$ induces an isomorphism on  $\mathcal K$.
Thus for  $Z\in \mathcal K$ with $Z_o\in \Delta_o$, by \eqref{decomposition}  there exist  $Y\in \g$ and $Z^\prime\in \mathcal K_o$ such that $f_* Z = Y^\dag + Z^\prime$. Since $f(o)=o$, we have that 
$$Y^\dag_o =(f_* Z)_o= df_o Z_o
= Z_o.
$$
So 
$$
f_*Z= Z+(Y^\dag -Z) + Z^\prime \in Z+\mathcal K_0.
$$
Hence (i) is proven.

Regarding the proof of (ii), we proceed by induction on $j$. The case $j=-1$ is given by (i).
Now suppose that (ii) holds for every index from $-1$ to $j-1$ and choose $Z\in \mathcal K_j$. Clearly $(f_*Z -Z)_{o}=0$. We are left to prove that $[f_*Z -Z,Y^\dag] \in \mathcal K_j$, for every $Y\in V$. Using (i) and the induction hypothesis, we have
\begin{align}\label{bothhandsnext}
[f_*Z-Z,Y^\dag]&= [f_*Z,Y^\dag]-[Z,Y^\dag]\\
&\in [f_*Z,f_*Y^\dag +\mathcal K_0] - [Z,Y^\dag] \nonumber\\
&= f_*[Z,Y^\dag]+[f_*Z,\mathcal K_0]-[Z,Y^\dag] \subseteq \mathcal K_j + [f_*Z,\mathcal K_0]\nonumber
\end{align}
Notice that if  $f_*Z\in \mathcal K_j$, then  Lemma \ref{gradationofkilling} together  with \eqref{bothhandsnext} imply (ii).
In order to prove that  $f_*Z\in \mathcal K_j$, we show that if $f_*Z\in \mathcal K_l$ for some $l <j$, then $f_* Z\in \mathcal K_{l+1}$. First, if $f_* Z\in \mathcal K_l$, then $[f_*Z,\mathcal K_0]\subseteq \mathcal K_l$. Then \eqref{bothhandsnext} implies 
$$
[f_*Z-Z,Y^\dag]\in \mathcal K_j + \mathcal K_l\subseteq \mathcal K_l.
$$
Therefore $f_*Z-Z\in \mathcal K_{l+1}$, whence 
$$
f_*Z\in Z+\mathcal K_{l+1}\subseteq \mathcal K_j +\mathcal K_{l+1}\subseteq \mathcal K_{l+1},
$$
because $j>l$. This concludes the proof.
\qed

\subsection{Unique extension for Killing vector fields}\label{tanaka}
In this section we shall prove  that an isometry defined among some open subsets of a subRiemannian Carnot group induces an isomorphism of the space $\mathcal K$ of Killing vector fields.  
This fact  will be important in order  to apply the results of the previous section to the proof of Theorem \ref{localisometries}.
  A key tool throughout this section is Tanaka's prolongation theory. On the one hand, we shall use a finiteness criterium \cite[Corollary 2, page 76]{Tanaka} to show finite dimensionality of Tanaka's prolongation, which will be in turn isomorphic to $\mathcal K$. On the other hand,   the discussion in \cite[Section 6]{Tanaka} will imply that a locally defined Killing vector field is uniquely extended to an element of $\mathcal K$.

We only recall what of Tanaka's theory is essential for our purposes. For an insight, we refer the reader  to  \cite{Tanaka, Yamaguchi}. A recent survey can also be found in \cite{Ottazzi_Warhurst}.

A linear self map $u:\g\to\g$ of a Lie algebra $\g$ is a derivation if 
$$u[X,Y]=[u(X), Y]+[X,u(Y)], \qquad \forall X,Y\in \g.$$
Given a nilpotent Lie algebra $\g$ with stratification $V_1\oplus\ldots\oplus V_s$, we set
$${\rm Der}_0(\g):=\{u:\g\to\g\text{ derivation  }\mid  u(V_j)\subseteq V_j, \forall j=1,\ldots, s\}.$$
Once we fix a scalar product on $V_1$, we denote by $O(V_1)$ the Lie group of linear isometries of $V_1$, and by $\mathfrak o (V_1)$ its Lie algebra.
In particular, $\mathfrak o (V_1)\subseteq \mathfrak {gl} (V_1)$.

We recall now the definition of the   Tanaka prolongation of the Lie algebra $\g$ with respect to the subspace
$$\g_0:= \{u\in {\rm Der}_0(\g) \mid u|_{V_1}\in \mathfrak o (V_1)\}.$$
Set $\g_j:=V_{-j}$, for $j\in \{-s,\ldots, -1\}$.
By induction, if $k\geq 1$, define $$\g_k:=\{u\in \bigoplus_{j=-s}^{-1} \g_{j+k}\otimes \g_j^*\mid u[X,Y]=[u(X), Y]+[X,u(Y)], \forall X,Y\in \g\},$$
where, if $X\in \g$ and $u\in \g_k$ with $k\geq 0$, we set
\begin{equation}\label{bracket_ux}
[u,X]:=u(X).
\end{equation}
The {\em Tanaka prolongation} of the Lie algebra $\g$ with respect to $\g_0$ is the sum
$${\rm Prol}(\g):=  \bigoplus_{j\geq-s} \g_{j}.$$
The Lie algebra structure of $\g$ is  extended to   ${\rm Prol}(\g)$  by \eqref{bracket_ux} and by the inductively defined formula
\begin{equation}\label{bracket_uvx}
[u,v](X):=[u(X),v]+[u,v(X)], \quad \forall u\in \g_i, v\in\g_j,  X\in \g.
\end{equation}
\begin{lemma}\label{finiteness}
${\rm Prol}(\g)$ is a graded finite dimensional Lie algebra. Namely
$[\g_i,\g_j]\subseteq \g_{i+j}$, $\forall i, j\geq -s$.  \end{lemma}
\proof
The grading property and the Jacobi identity easily follow  from the definition of the Lie bracket.
By \cite[Corollary 2, pag 76]{Tanaka},
for showing finite dimensionality it is enough to prove that the space
$$\mathfrak h_0:= \{u\in \g_0 \mid u[X,Y]=0, \forall X,Y\in \g\}$$
has finite prolongation in the sense of Singer and Sternberg, see \cite{Singer-Sternberg} or \cite{Kobayashi}.
The set $\mathfrak h_0$ can be identified with a subalgebra of $\mathfrak o (V_1)$.
By an easy argument  the first prolongation of $\mathfrak o (V_1)$ is trivial, 
see \cite[page 8]{Kobayashi}.
\qed

We observe that the definition of prolongation that we provided above adapts to every choice of a subalgebra
of ${\rm Der}_0(\g)$.
Tanaka constructed these prolongation algebras to describe the infinitesimal generators
of mappings preserving some nonintegrable geometric structures. 
If, e.g., we prolong with respect to all of ${\rm Der}_0(\g)$, then we would obtain a characterization of those vector fields generating 
flows of contact mappings. More precisely, one has that the prolongation algebra and the corresponding space of vector fields
identify as Lie algebras whenever the prolongation if finite.

In our case, since  ${\rm Prol}(\g)$ is finite by Lemma \ref{finiteness}, Tanaka's method leads to
the identification of Lie algebras
\begin{equation}\label{tanaka-killing}
{\rm Prol}{(\g)} \cong \mathcal K.
\end{equation}
Now we observe that this same identification also applies to vector fields that are defined only locally.
\begin{definition} (Local Killing vector fields $\mathcal K_\Omega$)
Let $(\G,d)$ be a Carnot group. Let $Z$ be a vector field on $\Omega\subseteq \G$.
We say that $Z$ is a {\em Killing vector field on} $\Omega$ if there exist $t_0>0$ and $\Omega_0\subset \Omega$ open set such that, for all $t\in [0,t_0]$, the flow $\Phi_Z^t|_{\Omega_0}$ at time $t$ is an isometric embedding of $\Omega_0$ inside $\G$, with respect to the Carnot-Carath\'eodory distance on $\G$.
We denote by $\mathcal K_\Omega$ the collection of all  Killing vector fields on an open set $\Omega$ of $\G$.
\end{definition}
The idea behind the isomorphism of $\mathcal K_\Omega$ and $\mathcal K$ is the following.
To any vector field  in $\mathcal K_\Omega$ 
we  associate an element of ${\rm Prol}(\g)$, the Tanaka prolongation of $\g$ with respect to $\g_0$.
By Lemma \ref{finiteness}, ${\rm Prol}(\g)$ is a finite dimensional Lie algebra.
One should think that the elements in ${\rm Prol}(\g)$ are the coefficients of the Taylor expansion at a fixed point of the coordinates of the vector field.
A vector field $Z$ is in $\mathcal K_\Omega$ if and only if the coefficients of $Z$ with respect to a fixed basis of left-invariant vector fields satisfy a particular system of PDE's, whose
 solutions are polynomials.
 Hence the coefficients of $Z$ are polynomials.
The same polynomials define an extension  of $Z$ for which we abuse the notation $Z$.
Since the system of PDE's is in fact polynomial and linear, then the value of this operator on the coefficients of $Z$ is a polynomial that is null on an open set. Therefore it is $0$ everywhere.
Hence, the vector field $Z$ is still a solution.
Thus, $Z\in \mathcal K$.
\begin{lemma}\label{lemma_A5} Let $\G$ be a subRiemannian Carnot group.
If $\Omega$ is a connected open subset of $\G$, then the restriction function from $\mathcal K=\mathcal K_\G$ to $\mathcal K_\Omega$ is a bijection.
\end{lemma}

\proof
One of the fundamental point in Tanaka's work is that a locally defined infinitesimal generator of mappings preserving some geometric structure can be extended to the whole group, whenever the prolongation describing 
that type of mappings is finite.
From Section 6 in \cite{Tanaka}, it follows in particular that to any vector field  in $\mathcal K_\Omega$ we can associate an element of ${\rm Prol}(\g)$, the Tanaka prolongation of $\g$ with respect to $\g_0$. So  \eqref{tanaka-killing} concludes the proof.
\qed 
\begin{remark}
Let $f$ be an isometry between two connected open sets of a Carnot group. By Theorem \ref{capcow},  the map $f$ is smooth. By Lemma \ref{lemma_A5},
the vector field  $f_*Z$, which is in  $\mathcal K_\Omega$,  uniquely extends to a vector field in  $\mathcal K$. Therefore, the isometry $f$ induces an isomorphism of $\mathcal K$.
\end{remark}

\section{Isometries of subFinsler homogeneous spaces}

This section is devoted to the proof of Theorem \ref{globalisometries} and Theorem \ref{localisometries}. 
We prove the theorem on global isometries immediately in Section \ref{globalthm}. 
For the proof of Theorem \ref{localisometries}, we first need 
to consider the situation of a distance  defined by a scalar product on $V_1$ and isometries defined on connected subsets (see Theorem \ref{local-connected}). Then we consider the general case in Section \ref{endproof}.

\subsection{Global isometries of subFinsler homogeneous spaces}\label{globalthm}
We prove now Theorem \ref{globalisometries}. 
We shall only make use of 
Corollary \ref{compactness} and Lemma \ref{gradinggeneral} of the previous discussion.

 \proof[Proof of Theorem \ref{globalisometries}]
From Corollary \ref{compactness}, it follows that $f$ and $h$ are smooth.
%Let $\phi$ be the blow up of $f$ at $e$, which exists since $f\in C^\infty$ and contact. 
%By \cite{Margulis-Mostow}, the map $\phi$ is a group isomorphism of the tangent cone and, moreover, it is an isometry, being the limit of isometries.
%In particular we have that 
%$(dh)_e|_{V}=d\phi_e |_{V}=df_e |_{V}.$
We  suppose that  we have $f(o)=h(o)$ and $(dh)_o|_{\Delta_o}=df_o |_{\Delta_o}.$
We plan to show that $f=h$.

Up to replacing $f$ with $h^{-1}\circ f$, we may assume that  $f(o)=o$ and that $df_o$ is the identity on $\Delta_o$.
%\begin{equation*}%\label{ident_on_V1}
%df_oY=Y,\qquad \forall Y\in \Delta_o.
%\end{equation*}
Let $f_*$ be the push forward operator on vector fields.
Since $f$ is smooth and globally defined, the map $f_*$ is a Lie algebra isomorphism on  $\mathcal K$, the space of  Killing vector fields.
We shall prove that $f_*Y^\dag=Y^\dag$, for all $Y\in \g$.
First, pick $Y\in V$ and 
%By Lemma \cite{lemma_B1}, there exists $R_X\in C_0$ such that $\alpha(X)=X+R_X.$
assume by contradiction that $f_*Y^\dag\neq Y^\dag$. Then, using Lemma \ref{gradinggeneral}, there exists $Z\in \mathcal K_j$, with $j\geq 0$ and $Z\neq 0$, such that
$$f_*Y^\dag\in Y^\dag+Z + \mathcal K_{j+1}.$$
Using again  Lemma \ref{gradinggeneral}, we have
 \begin{eqnarray*}
f_*^2 Y^\dag &=& f_*f_*Y^\dag \\
&\in& f_*Y^\dag +f_*Z + f_*\mathcal K_{j+1}\\
&\subseteq&  Y^\dag+Z + \mathcal K_{j+1}+Z + \mathcal K_{j+1} + \mathcal K_{j+1}= Y^\dag+2Z + \mathcal K_{j+1}.
 \end{eqnarray*}
By iteration, we get
\begin{equation}\label{iteration}
f_*^n Y^\dag  \in  Y^\dag+n Z + \mathcal K_{j+1}, \qquad \forall n\in \N.
\end{equation}
Notice that, for all $K\in \mathcal K$, 
$$f_*^nK\in \{g_*K\mid g\in {\rm Iso}_o(M)\}.$$
Since ${\rm Iso}_o(G)$ is compact (see Corollary \ref{compactness}), the family $\{g_*\}_{g\in {\rm Iso}_o(G)}$ is a collection of bounded operators.
So, on the one hand, $f_*^n Y^\dag  $ belongs to a bounded set.
On the other hand, any bounded set has empty intersection with the affine  space $Y^\dag+nZ + \mathcal K_{j+1}$, for $n$ big enough.
Hence, \eqref{iteration} is contradicted.
Thus $f_*Y^\dag=Y^\dag$ for every $Y\in V$.
Notice that formula \eqref{bracket_dag} implies that vector fields of the form $Y^\dag$ as $Y$ varies in $V$ bracket-generate all vector fields $Y^\dag$ with $Y\in \g$. 
Since $f_*$ commutes with the bracket of vector fields, we conclude that 
$$
f_*Y^\dag=Y^\dag,$$
for all $Y\in V$.
This implies
$$
f(\exp (tY)\cdot f^{-1}(gH))= \exp (tY) gH
$$
for every $Y\in \g$. In particular, choosing $g=o$ and since 
 $f(o)=o$, we obtain
 \begin{equation}\label{identity}
 f(\exp (tY) H)=\exp (tY) H,
 \end{equation}
 for every $Y\in \g$.
 Let now $U\subseteq G$ be a neighborhood of $e$ with the property that for every $g\in U$ there exists $Y\in \g$ such that $\exp Y =g$. Then \eqref{identity} implies that $f$ is the identity when restricted to the left cosets of $U$. Since $M$ is supposed to be connected, we conclude that $f$ is the identity everywhere.
 \qed

\subsection{Isometries of open sets in Carnot groups}\label{localthm}

In this section we prove Theorem \ref{localisometries} in a particular case. % that the distance is defined by a scalar product on $V_1$. In other words, 
Namely, we show that any isometry defined between two connected and open subsets of a Carnot group $\G$ endowed with a subRiemannian metric is affine. 
%Namely, it is the composition of a translation with a group isomorphism of $\G$.

The structure of the proof of the following theorem is similar to that one of Theorem \ref{globalisometries}. 
Since now we deal with isometries defined on subsets, we shall rely on  \cite{Capogna-Cowling} for the smoothness of such maps. Moreover, Lemma \ref{lemma_A5} and the observation thereafter allow us to extend
any local Killing vector field to a global one. Consequently an isometry defined on an open subset of a Carnot group $\G$  will induce an isomorphism of $\mathcal K$.

\begin{theorem}\label{local-connected}
Let $(\G,d)$ be a subRiemannian Carnot group.
Let $\Omega_1,\Omega_2\subset \G$ be two {\em connected} open sets.
Let $f:\Omega_1\to \Omega_2$ be an isometry.
If $f(e)=e$, then $f$ is the restriction to $\Omega_1$ of a group isomorphism of $\G$.
\end{theorem}

\proof
We have that  $f$ is analytic from Theorem \ref{capcow}. Let $\phi$ be the blow up of $f$ at $e$, i.e., the Pansu differential at the identity, see
 \cite{Warhurst-2}.
Notice that 
$d\phi_e |_{V_1}=df_e |_{V_1}.$
By \cite{Pansu}, the map $\phi$ is a group isomorphism and, moreover, it is an isometry, being the limit of isometries.
We plan to show that $f=\phi|_{{\Omega_1}}$.
Up to replacing $f$ with $\phi^{-1}\circ f$, we may assume that 
\begin{equation*}%\label{ident_on_V1 bis}
df_eY=Y,\qquad \forall Y\in V_1.
\end{equation*}
We now prove that $f_*Y^\dag=Y^\dag$, for all $Y\in \g$.
Let $Y\in V_1$. 
Assume by contradiction that $f_*Y^\dag\neq Y^\dag$.
By means of Lemma \ref{lemma_A5} and the remark thereafter, the map $f_*$ induces a Lie algebra isomorphism of the Killing vector fields $\mathcal K$. Therefore, proceeding as in the first part of the proof of Theorem \ref{globalisometries}, we  show that 
\begin{equation}\label{iterationbis}
f_*^nY^\dag  \in  Y^\dag+nZ + \mathcal K_{j+1}, \qquad \forall n\in \N,
\end{equation}
for  some nonzero element $Z\in \mathcal K_j$.
Let $\mathcal I$ be the group of isometries preserving a ball on which $f$ is defined.
Notice that, for all $K\in \mathcal K$, 
$$f_*^nK\in \{g_*K\mid g\in \mathcal I\}.$$

By (iii) of Corollary \ref{Hilbert} the differentials $g_*$ depend smoothly on $g\in \mathcal I$. Since $\mathcal I$ is compact, then the family $\{g_*\}_{g\in \mathcal I}$ is a collection of bounded operators.
In particular, the set
$\{g_*Y^\dag\mid g\in \mathcal I\}$ is bounded.
So, on the one hand, $f_*^n Y^\dag  $ belongs to a bounded set.
On the other hand, any bounded set has empty intersection with the affine  space $Y^\dag+nZ + \mathcal K_{j+1}$, for $n$ big enough.
Hence, \eqref{iterationbis} is contradicted.
Thus $f_*Y^\dag=Y^\dag$, for every $Y\in V_1$ . Since $f_*$ is a Lie algebra automorphism, we conclude that $f_*$ is the identity on right-invariant vector fields.  But $f(e)=e$, so $f$ is the identity.
 \qed

\subsection{Generalization to nonconnected set and to subFinsler Carnot groups}\label{endproof}
We shall extend Theorem \ref{local-connected} to the case of subsets of $\G$ that are not necessarily connected.
In order to do that, we use the following theorem.
\begin{theorem}[{\cite[Theorem 1]{Agrachev09}}]\label{agra}
Let $M$ be a (analytic) subRiemannian  manifold and set $q_0\in M$. Then there exists an open and dense subset $\Sigma_{q_0}\subseteq M$ such that for any $q\in \Sigma_{q_0}$ there exists a unique length minimizing curve $\gamma$ connecting $q_0$ to $q$. Moreover, such $\gamma$ is analytic.
\end{theorem}
The following result holds for general subRiemannian manifolds. We show that an isometry is completely determined on its behavior on an open set.
\begin{proposition}\label{agraconsequence}
Let $M$ be a (analytic) subRiemannian  manifold. Let $f:{\Omega_1}\to \Omega_2$ be an isometry among two open sets in $M$. Assume that $f$ is the identity on an open subset of $\Omega_1$. Then $f$ is the identity.
\end{proposition}
\proof
%Let now $f$ be an isometry among two open subsets of $\G$. 
%By Theorem \ref{local-connected}, $f$ restricted to any connected component of its domain is a group isomorphism. 
%We show that $f$ is the same group isomorphism on every pair of connected components $\Omega_1$ and $\Omega_2$. Without loss of generality, we may assume that $f|_{\Omega_1}$ is the identity. 

Let $\Omega\subseteq \Omega_1$ be the open subset such that $f|_{\Omega}$ is the identity.
Pick $q\in \Omega_1$. According to the notation in Theorem \ref{agra}, consider $\Sigma=\Sigma_q \cap \Sigma_{f(q)}$. Fix $p\in \Sigma\cap \Omega$. Since $p\in\Sigma_q$, Theorem \ref{agra} implies that there exists a unique and analytic length minimizing curve $\gamma$ such that $\gamma(0)=p$ and $\gamma(1)=q$. Since $\Omega$ is open, one can choose  $s_0\in (0,1)$ such that  $p^\prime:=\gamma(s_0)\neq p$ and $\gamma(s_0)\in \Omega$. Denote by $\rho$ a length minimizing curve such that $\rho(0)=p^\prime$ and $\rho(1)=f(q)$. Let $\tilde \gamma$ be the curve formed by joining $\gamma|_{[0,s_0]}$ with $\rho$. We claim that $\tilde \gamma$ minimizes the length between $p$ and $f(q)$. %Indeed, notice that $p,p^\prime$ and $q$ lie on geodesics.
Indeed, since $f(p)=p$, $f(p^\prime)=p^\prime$, and $f$ is an isometry, we have
\begin{align*}
d(p,f(q))&= d(p,q)= d(p,p^\prime)+ d(p^\prime,q)\\
&= d(p,p^\prime) + d(p^\prime,f(q)).
\end{align*}
So $\tilde \gamma$ realizes the distance from $p$ to $f(q)$.
Since $p\in \Sigma_{f(q)}$, it follows that $\tilde \gamma$ is analytic. Since $\gamma$ and $\tilde \gamma$ coincide on an interval, they are both analytic, and have the same length, they coincide.  In particular, $q=f(q)$.
\qed

Using the proposition above, together with Theorem \ref{local-connected},  we obtain Theorem \ref{localisometries} in the case of subRiemannian Carnot groups. 
At this point,   Theorem \ref{localisometries}  will be proved once we extend the results to the case of a subFinsler metric. This is the content of the following lemma.
\begin{lemma}\label{subfinsler}
Let $f$ be a $C^1$ isometry among open subsets of $\G$ with respect to a left-invariant subFinsler distance $d_{SF}$. Then there exists a left-invariant subRiemannian distance $d_{SR}$ with same horizontal bundle as $d_{SF}$ such that $f$ is an isometry with respect to $d_{SR}$.
\end{lemma}
\proof

The proof is an application of John's Ellipsoid Theorem, see \cite{John_ellipsoid}. Let $\tilde V_p$ be the horizontal bundle at $p$. Denote by $K_p=\{v\in \tilde V_p\mid \| v\|\leq 1\}$, where $\|\,\cdot\,\|$ is the norm defining $d_{SF}$.
John's Ellipsoid Theorem states that there exists a unique maximal ellipsoid $E_p$ contained in $K_p$.

 Let $f$ be any $C^1$ $d_{SF}$-isometry. We claim that for any $p$ in the domain of $f$, we have
 \begin{equation}\label{isoellipse}
 df_p(E_p)=E_{f(p)}.
 \end{equation}
Indeed, $df_p$ restricts to a linear isometry between $(\tilde V_p,\|\,\cdot\,\|)$ and $(\tilde V_{f(p)}, \|\,\cdot\,\|)$. In particular, $df_p(K_p)=K_{f(p)}$ and $df_p(E_p)$ is an ellipsoid contained in $K_{f(p)}$.
Since $E_{f(p)}$ is the unique maximal ellipsoid, it follows that $df_p(E_p)\subseteq E_{f(p)}$. Since $df_p^{-1}$ also restricts to a linear isometry, we obtain the reverse inclusion of \eqref{isoellipse}. In particular, choosing $f$ to be a left translation $L_p$, with $p\in \G$, we have that
$
E_p= (dL_p)_e E_e
$.
Therefore, $\{E_p\}_{p\in\G}$ define a left-invariant scalar product on $\tilde V_p$ which in turn gives a left-invariant subRiemannian distance on $\G$. Moreover, the equation \eqref{isoellipse} implies that any $C^1$ isometry with respect to  $d_{SF}$ is also an isometry with respect to $d_{SR}$.
\qed

%We now have all the tools to prove Theorems \ref{localisometries} and \ref{globalisometries}.
\subsubsection{Proof of Theorem \ref{localisometries}}
Let $(\G,d_{SF})$ be a subFinsler Carnot group. Let $\Omega_1,\Omega_2\subseteq \G$ be two open sets and $f:\Omega_1\to \Omega_2$ an isometry. Arguing as in the proof of Lemma \ref{subfinsler}, we apply John's Ellipsoid Theorem and obtain a subRiemannian distance $d_{SR}$ on $\G$ satisfying the following properties:
 $(\G,d_{SR})$ is a subRiemannian Carnot group and any $C^1$ isometry of $\G$ with respect to $d_{SF}$ is an isometry with respect to $d_{SR}$.

Since $d_{SR}$ is biLipschitz equivalent to $d_{SF}$, it follows that $f$ is biLipshitz with respect to $d_{SR}$. By Pansu Theorem \cite{Pansu}, the blow up exists a.e. and it is a group isomorphism, hence $C^1$. Since any blow up of $f$ is an isometry of $(\G,d_{SF})$, Lemma \ref{subfinsler} implies that the blow up of $f$ at a.e. point is an isometry of $(\G,d_{SR})$. By Theorem \ref{capcow}.(ii), the map $f$ is analytic. Again using Lemma  \ref{subfinsler}, we have that the map $f$ is an isometry with respect to $d_{SR}$. Up to composing with a translation, we may assume $f(e)=e$. 
By Theorem \ref{local-connected}, we obtain that on the connected component $\Omega$ of $\Omega_1$ containing $e$, the map $f$ is a group isomorphism $\phi$. Then the map $\phi^{-1}\circ f$ is an isometry that is the identity on $\Omega$. By Proposition \ref{agraconsequence}, we get that $\phi^{-1}\circ f$ is the identity on $\Omega_1$, which finishes the proof.
\qed

%%%%%%%%%%%%%%%%%%%%%%%%

\subsection{Afterthoughts}\label{afterthoughts}
Not even in Euclidean space it is true that all group isomorphisms are isometries. 
However, an automorphism of a subFinsler Carnot group is an isometry if and only if its differential  preserves the first stratum (ad hence all strata) and restricted to the first stratum preserves the norm defining the subFinsler distance.
Hence we have a complete description of local isometries of subFinsler Carnot groups.

Unlikely in the Euclidean space, Theorem \ref{localisometries} cannot be generalized to arbitrary subsets. Here we present a counterexample.  We take the subRiemannian Heisenberg group $(\mathbb H, d_{SR})$ and we define exponential coordinates $(x,y,z)$ with respect to the basis of its Lie algebra given by vectors $X, Y$ and $Z$. The only nonzero bracket relation is $[X,Y]=Z$.  Consider the three coordinate axes,
namely,
$$E:=\exp (\R X) \cup \exp (\R Y)  \cup \exp (\R Z).$$
Then the map 
 $(x,y,z)\mapsto (x,y,-z)$
%$(x,y,z)\mapsto (x,y,-z)$
%
%
%that is the identity on $\exp (\R X) \cup \exp (\R Y)$ and sends $\exp (t Z)$ to $\exp (-t Z)$
 is an isometry of $E$ into itself. However, this map is not the restriction of a group isomorphism.
 
  Given an isometry $f$ of a subFinsler homogeneous space $G/H$, the differential of the blow-up at a point $p$ equals the differential 
 at $p$, when they are both restricted to $\Delta_p$. Therefore
 Theorem \ref{globalisometries}  claims that ${\rm Iso}_o(G/H)$  injects into ${\rm Iso}_o((G/H)_o)$, where $(G/H)_o$ denotes the Gromov tangent cone of $G/H$ at  $o$, which is a Carnot group. However, it is not true that isometries of $(G/H)_o$ are always  blow-ups of isometries of $G/H$. In fact, we can find counterexamples even in the domain of Riemannian Lie groups. Take for instance the three dimensional Heisenberg group, endowed with a Riemannian distance. We denote it by $(\mathbb H,d_{R})$. 
 Then its tangent cone at every point is the Euclidean $3$-space, 
 which contains all the rotations  among its isometries. However,  rotations with respect to horizontal lines are not  isometries for $(\mathbb H,d_{R})$. This follows from the observation that  ${\rm Iso}(\mathbb H,d_{R})={\rm Iso}(\mathbb H, d_{SR})$. The identification of the two isometry groups 
rests upon the study of length minimizing curves: for both metric models of $\mathbb H$ the only infinite geodesics are the $1$-parameter groups corresponding to horizontal vectors. Since isometries must preserve infinite geodesics, it  follows that the horizontal space is preserved by the differential of any isometry.

We notice that the statements of
Theorem \ref{localisometries} and Theorem \ref{globalisometries}  become equivalent if $G/H=\G$ and if $\Omega_1= \G$. If this is not the case, we cannot conclude that a global isometry of a subFinsler Lie group $G$ is affine. 
 Indeed, take  the three dimensional sphere $\mathbb S^3$ viewed as the space of quaternions with euclidean norm equal to one. The manifold $\mathbb S^3$ is then a Riemannian Lie group. It is easy to check that the inversion map $p\mapsto p^{-1}$ on $\mathbb S^3$ is an isometry that is not a group isomorphism.

For general Lie groups, isometries between open sets might not be 
restrictions of global isometries. For example,  
every point 
in the flat cylinder $\R\times \mathbb S^1$
has a neighborhood isometric to a disk in $\R^2$. Hence, all rotations are isometries of such a neighborhood. Of course, not all of them extend to global isometries.

Last, the fact that isometries of a Carnot group $\G$ are affine maps implies that $\mathcal K_j =\{0\}$ for every $j\geq 1$. As a by-product, the Tanaka prolongation ${\rm Prol}(\g)$ defined in Section \ref{tanaka} is $\g \oplus \g_0$. Although Tanaka's method proves the finiteness of ${\rm Prol}(\g)$ (see Lemma \ref{finiteness}) , we are not aware of a direct method to show that in general ${\rm Prol}(\g)=\g \oplus \g_0$.

 \bibliography{general_bibliography}
\bibliographystyle{amsalpha}

\end{document}